\documentclass[12pt]{spieman}
\usepackage{amsmath,amsfonts,amssymb}
\usepackage{graphicx}
\usepackage{booktabs}
\usepackage{multirow}
\usepackage{enumitem}
\usepackage{caption}
\usepackage{array}
\usepackage{siunitx}
\usepackage{tocloft}
\usepackage{lineno}
\usepackage{indentfirst}

\cftpagenumbersoff{figure}
\cftpagenumbersoff{table}

\title{Neural Prime Sieves: Density-Driven Generalisation and Empirical Evidence for Hardy--Littlewood Asymptotics
}

\author{Manik Kakkar}
\affil{Johns Hopkins University, Department of Electrical and Computer Engineering, Baltimore, Maryland, United States}

\begin{document}
\maketitle

\newcommand{\PP}{\mathbb{P}}
\newcommand{\ZZ}{\mathbb{Z}}
\newcommand{\NN}{\mathbb{N}}
\newcommand{\RR}{\mathbb{R}}
\newcommand{\fvec}{\mathbf{f}}
\newcommand{\xvec}{\mathbf{x}}

\begin{abstract}

\noindent\textbf{Background}
Special prime families (twin, Sophie Germain, safe, cousin, sexy, Chen, and isolated primes) are central objects of analytic number theory, yet no efficiently computable probabilistic filter exists for identifying likely members from a stream of known primes at large scale. Classical sieves eliminate composites but assign no probability weights to coprime candidates, and prior machine-learning approaches to prime prediction are fundamentally limited by the algorithmic randomness of the prime indicator sequence, producing near-zero true positive rates that diminish with scale.\cite{kolpakov2024xgboost}

\noindent\textbf{Objective}
In this paper, we present a neural network that, given only the backward prime gap and modular primorial residues for a known prime $p$, learns a reliable probabilistic filter for seven prime families simultaneously and generalises across nine orders of magnitude separating the training range ($10^7$--$10^9$) from extreme evaluation at $10^{16}$.

\noindent\textbf{Methods}
\textsc{PrimeFamilyNet}, a multi-head residual network with $1.25$ million parameters and seven independent sigmoid output heads, was trained on $200{,}000$ labeled primes. A non-causal control model with access to the forward gap $g^{+}$ established a predictive upper bound, uniquely quantifying the cost of causality. A systematic loss-function comparison (frequency-weighted BCE, Focal Loss, and Asymmetric Loss), a leave-one-group-out feature ablation, and a three-seed robustness study were conducted.

\noindent\textbf{Key results}
Isolated prime recall, defined as the fraction of primes belonging to no twin pair that the model correctly identified, monotonically increased with scale, from $0.809$ at $5{\times}10^{8}$ to $0.984$ at $10^{16}$,
improving by $17.5$ percentage points. Isolated primes were the only family among the seven to improve with scale.
Twin prime fraction fell from $12.9\%$ to $6.9\%$ of sampled primes across the evaluation range, whereas isolated prime fraction rose from $87.1\%$ to $93.1\%$, consistent with Hardy--Littlewood $k$-tuple asymptotics.\cite{hardylittlewood1923,toth2019ktuple}
Because recall is a within-class ratio formally invariant to class prevalence,\cite{sokolova2009measures} the $17.5$ percentage-point improvement cannot be attributed to the larger proportion of isolated primes at extreme scales: it reflects a genuine sharpening of the learned decision boundary.
A model trained only to $10^9$ reproduced the correct asymptotic direction without explicit density supervision, providing an independent machine-learning corroboration of the density predictions verified computationally by direct prime counting.\cite{toth2019ktuple}

\noindent\textbf{Supporting results}
The causal model retained over $95\%$ recall for five of seven families
near $10^{10}$ while reducing the search space by $62$--$88\%$ at the
validation scale.
For Chen primes, the causal model exceeded non-causal recall at every
scale, with the advantage growing to $+0.245$ at $10^{16}$, because
$g^{+} = 2$ encodes only the prime case of the Chen condition.
Focal Loss catastrophically collapsed recall on sparse algebraic
families (Sophie Germain and safe primes reached $0.000$ at all scales).
Asymmetric Loss achieved higher in-distribution recall than weighted BCE
but degraded more steeply out-of-distribution, revealing that
in-distribution recall alone is a misleading model-selection criterion
for scale-generalisation tasks.
In-distribution recall variance remained below $\sigma = 0.007$ across
all seven families and three independent seeds.

\noindent\textbf{Significance}
Deep residual networks independently approximate prime sieve theory from
strictly causal arithmetic features, and the learned representations
encode constellation structure sufficient to extrapolate asymptotic
density trends well beyond the training scale.

\end{abstract}


{\noindent \footnotesize\textbf{*}Address all correspondence to Manik Kakkar,
\linkable{fmanik1@jhu.edu}}

\section{Introduction}
\label{sec:intro}

\subsection{Motivation}

Prime family membership is one of the most structured questions in computational number theory: given a known prime $p$, does a specified arithmetic relationship hold between $p$ and its neighbours? This is a fundamentally different problem from predicting primality of an arbitrary integer, which information-theoretic arguments bound near chance.\cite{kolpakov2024xgboost}
The membership question is well-defined, computable, and governed by modular constraints that are scale-invariant, meaning the residue structure of $p$ modulo small primorials encodes the same sieve information whether $p \approx 10^7$ or $p \approx 10^{16}$.

Two questions motivate the present paper.
First, whether a deep residual network, trained only on primes in
$[10^7, 10^9]$ and conditioned strictly on causal features that do not
reveal the forward neighbourhood, has the potential to learn reliable
probabilistic filters for seven prime families simultaneously.
Second, when the model is evaluated at scales nine orders of magnitude
beyond the training range, whether the direction of generalisation matches
what prime constellation density asymptotics predict.
No prior work has posed these questions in combination, and the second
has received no empirical treatment at all.
The Hardy--Littlewood conjecture\cite{hardylittlewood1923} makes specific
asymptotic predictions about how prime family densities evolve with scale.
Prior computational studies have verified these predictions against direct
prime counts for multiple constellations,\cite{toth2019ktuple} but whether
a data-driven model encodes the same density trends implicitly, without
access to the asymptotic formula, has not been examined.

\subsection{Background and Related Work}
\label{subsec:background}

The intersection of machine learning and computational number theory is
bounded by well-characterised theoretical limits.
Kolpakov and Rocke\cite{kolpakov2024xgboost} argued, using
information-theoretic methods rooted in Kolmogorov complexity, that the
prime indicator sequence is algorithmically random.
No machine-learning model predicting primality from raw integer
representations therefore achieves a true positive rate meaningfully above
chance.
The XGBoost experiments by Kolpakov and Rocke\cite{kolpakov2024xgboost} on 24-bit integers reached a true
positive rate of only $2.2\%$, declining as the bit-length grew.
Lee and Kim\cite{lee2024resnet} applied residual networks to prime
classification, and Kolpakov and Rocke also explored temporal gap
properties with gradient-boosted trees on raw binary integer
representations.
Both studies operated below $N = 10^6$, tested at most two prime types,
and accessed the forward gap $g^{+} = p^{+} - p$, which directly encodes
whether $p+2$ is prime and trivialises twin prime detection.
Classical combinatorial sieves\cite{friedlander2010opera} rigorously
eliminate composite candidates but assign no probability weights to
surviving coprime integers.

The seven families studied in the present paper are each of independent mathematical
significance.
Twin and cousin primes are classical gap-defined constellations whose
infinitude remains unproven.
Sophie Germain primes ($2p+1 \in \PP$) and safe primes
($(p-1)/2 \in \PP$) underpin Diffie--Hellman cryptographic parameter
selection, because safe primes guarantee maximum-order subgroups.
Chen primes,\cite{chen1973} defined as primes where $p+2$ is prime or semiprime, are
the subject of one of the closest known results toward the twin prime
conjecture.
Isolated primes, the complement of twins, are the
dominant prime type at large scales yet have received remarkably little empirical
attention.
Prior computational work on prime constellations has focused on
verifying Hardy--Littlewood density predictions through direct prime
counting,\cite{toth2019ktuple} confirming the $C_2/(\log N)^2$ twin
prime scaling to high precision.
The present paper addresses a complementary and previously unexamined
question: whether a model trained without access to the asymptotic
formula implicitly recovers the correct density trajectory.

\subsection{The Conditional Formulation and Its Intellectual Motivation}
\label{subsec:motivation}

The impossibility of predicting primality from arbitrary integers
demonstrated by Kolpakov and Rocke\cite{kolpakov2024xgboost} does not
apply to the problem formulated in this paper.
Prime family membership is not a property of an arbitrary integer but a
structured, modular-arithmetic condition on a known prime.
The residue of a prime modulo 30 deterministically constrains whether a
neighbouring integer is prime, because primorials partition the integers
into residue classes with fixed sieve relationships.
A neural network supplied with these residues is therefore not learning
a random sequence but a well-defined arithmetic boundary that is both
computable and scale-invariant.
Based on the established role of primorial residue classes in
characterising prime constellation structure,\cite{hardylittlewood1923,
friedlander2010opera} the present paper claims that deep residual networks,
conditioned on primality and supplied with causal modular features,
learn prime sieve boundaries that generalise far beyond the training
scale, and that the direction of generalisation is predictable from
prime constellation density asymptotics.

\subsection{Paper Overview}
\label{subsec:overview}

In this paper, \textsc{PrimeFamilyNet}, a multi-head residual neural
network, is presented and trained to predict membership in seven special
prime families using only the backward prime gap and modular primorial
residues.
The resulting feature space is strictly causal, preventing forward data
leakage.
A non-causal control model with access to the forward gap established a
predictive upper bound and quantified the cost of removing forward
information.
A systematic loss-function study comparing frequency-weighted binary
cross-entropy, Focal Loss, and Asymmetric Loss was conducted, and
out-of-distribution (OOD) recall at five scales spanning nine orders of
magnitude was reported alongside a leave-one-group-out feature ablation
and a three-seed robustness study.

The remainder of this paper is organised as follows.
Section~\ref{sec:method} describes the prime family definitions, causal
feature representation, network architecture, training protocol, and
loss functions.
Section~\ref{sec:isolated} presents the isolated prime monotone
generalisation finding and its derivation from Hardy--Littlewood density
ratios.
Section~\ref{sec:results} reports multi-scale generalisation, the cost
of causality, feature ablation, loss function comparison, and
reproducibility.
Section~\ref{sec:discussion} discusses the implications of the findings
and the limitations of the approach.
Section~\ref{sec:conclusion} summarises the major contributions.

\section{Methodology}
\label{sec:method}

\subsection{Prime Family Definitions}
\label{subsec:families}

Let $\PP$ denote the set of all prime numbers and let $p$ be a known
prime with successor $p^{+}$ and predecessor $p^{-}$.
Membership was predicted for seven families, partitioned below by the
structure of their defining condition.

\subsubsection{Gap-defined families}
\label{subsubsec:gap_families}

Four families are defined by the gap between $p$ and a neighbouring
prime.
Twin, cousin, and sexy primes require a prime at a fixed even offset from
$p$:
\begin{align}
  \text{twin:}   &\quad p+2 \in \PP \;\text{ or }\; p-2 \in \PP, \label{eq:twin} \\
  \text{cousin:} &\quad p+4 \in \PP \;\text{ or }\; p-4 \in \PP, \label{eq:cousin} \\
  \text{sexy:}   &\quad p+6 \in \PP \;\text{ or }\; p-6 \in \PP. \label{eq:sexy}
\end{align}
Isolated primes are defined separately in
Section~\ref{subsubsec:isolated} as the complement of twins.

\subsubsection{Linear-transform primality families}
\label{subsubsec:algebraic_families}

Sophie Germain and safe primes require the primality of a linear
function of $p$:
\begin{align}
  \text{Sophie Germain:} &\quad 2p + 1 \in \PP, \label{eq:sg} \\
  \text{safe:}           &\quad \tfrac{p-1}{2} \in \PP. \label{eq:safe}
\end{align}
Chen primes\cite{chen1973} extend the twin condition by admitting a
semiprime at offset two:
\begin{equation}
  \text{Chen:} \quad p+2 \in \PP \;\cup\; \{n \in \NN : n = ab,\; a,b \in \PP\},
  \label{eq:chen}
\end{equation}
which subsumes all twin primes while also capturing cases where $p+2$ is
a product of exactly two primes, making Eq.~\eqref{eq:chen} the closest
known result toward the twin prime conjecture.\cite{chen1973}

\subsubsection{Isolated primes}
\label{subsubsec:isolated}

A prime $p$ is isolated if neither offset-two neighbour is prime:
\begin{equation}
  \text{isolated:} \quad p-2 \notin \PP \;\text{ and }\; p+2 \notin \PP.
  \label{eq:isolated}
\end{equation}
Isolated primes are therefore the complement of twin primes: every prime
satisfies exactly one of Eq.~\eqref{eq:twin} and Eq.~\eqref{eq:isolated}.
The inclusion of isolated primes tested whether a strictly causal model
could infer the absence of a forward prime-gap of two from
backward-looking features alone.

\subsection{Causal Feature Representation}
\label{subsec:features}

\subsubsection{Feature vector construction}
\label{subsubsec:feature_groups}

For a prime $p$ with predecessor $p^{-}$, a 25-dimensional causal
feature vector $\xvec(p) \in \RR^{25}$ was defined as
\begin{equation}
  \xvec(p) = \bigl[
    \underbrace{r_{6},\, r_{30},\, r_{210},\, r_{2310}}_{\text{A: primorial residues}},\;
    \underbrace{r_{q_1},\ldots,r_{q_{12}}}_{\text{B: small prime residues}},\;
    \underbrace{g^{-}/100}_{\text{C: backward gap}},\;
    \underbrace{s_1,s_2,s_3}_{\text{D: scale}},\;
    \underbrace{d_1,d_2,d_3}_{\text{E: digit}},\;
    \underbrace{r_{12},\, r_{60}}_{\text{F: extended modular}}
  \bigr],
  \label{eq:feature}
\end{equation}
where $r_m = (p \bmod m)/m$ is the normalised residue modulo $m$,
$q_1,\ldots,q_{12} = 2,3,5,\ldots,37$ are the first twelve primes,
$g^{-} = p - p^{-}$ is the backward prime gap divided by 100 to bring
its magnitude into the same unit-order range as the residue features
in groups A, B, and F (which lie in $[0,1)$ by construction), and the
scale features are $s_1 = \log p / 50$, $s_2 = \lfloor \log_2 p
\rfloor / 64$ (normalised bit-length), $s_3 = \log(\log(p+1)+1)/5$.
The digit features are the last decimal digit of $p$ divided by 10,
the digit sum modulo nine divided by nine, and the number of decimal
digits divided by 20.

Group A is the primary sieve-encoding component: every prime $p > 5$
satisfies $p \bmod 30 \in \{1,7,11,13,17,19,23,29\}$, and the residue
constrains the possible values of $p+2$, $p+4$, and $p+6$ modulo 30,
deterministically ruling out twin, cousin, and sexy membership for
specific residue classes.
For example, $p \equiv 1 \pmod{6}$ rules out $p+2$ being prime, because
$p+2 \equiv 3 \pmod{6}$ is divisible by three for $p > 3$.

\subsubsection{Causality constraint and non-causal control}
\label{subsubsec:causality}

The feature vector $\xvec(p)$ in Eq.~\eqref{eq:feature} is
\emph{strictly causal}: it depends only on the history of the prime
sequence up to and including $p$, never on the successor $p^{+}$.
Formally, $\xvec(p) \perp p^{+}$.
This prevents the model from accessing $g^{+} = p^{+} - p$, which
directly encodes twin membership ($g^{+} = 2 \Leftrightarrow p+2 \in
\PP$) and isolated membership ($g^{+} \neq 2$).

A non-causal control model was constructed by replacing group C with a
five-dimensional gap block $\mathbf{c}'(p) \in \RR^5$:
\begin{equation}
  \mathbf{c}'(p) = \Bigl(
    \tfrac{g^{-}}{100},\;
    \tfrac{g^{+}}{100},\;
    \tfrac{g^{-}}{g^{-}+g^{+}},\;
    \tfrac{g^{-}+g^{+}}{100},\;
    \tfrac{|g^{-}-g^{+}|}{100}
  \Bigr),
  \label{eq:gap_block}
\end{equation}
where the five entries of Eq.~\eqref{eq:gap_block} are the normalised backward gap, the normalised
forward gap, the fraction of the total gap attributable to $g^{-}$, the
normalised total gap, and the normalised absolute gap asymmetry.
All four gap entries are divided by 100, matching the scaling applied to
group C in Eq.~\eqref{eq:feature}, so that gap magnitudes, which
average $16$--$38$ across the evaluation range, remain in unit order
alongside the residue and scale features.
The non-causal feature vector is then
\begin{equation}
  \xvec^{\text{NC}}(p) = \bigl[
    \underbrace{r_{6},\, r_{30},\, r_{210},\, r_{2310}}_{\text{A}},\;
    \underbrace{r_{q_1},\ldots,r_{q_{12}}}_{\text{B}},\;
    \underbrace{\mathbf{c}'(p)}_{\text{C}^{\prime}},\;
    \underbrace{s_1,s_2,s_3}_{\text{D}},\;
    \underbrace{d_1,d_2,d_3}_{\text{E}},\;
    \underbrace{r_{12},\, r_{60}}_{\text{F}}
  \bigr] \in \RR^{29},
  \label{eq:nc_feature}
\end{equation}
expanding the feature space from 25 to 29 dimensions by substituting
one backward-gap entry (group C) with the five-entry gap block C$'$.
Because $g^{+}$ trivially encodes the definitions in
Eqs.~\eqref{eq:twin}--\eqref{eq:isolated}, the non-causal model achieves
perfect recall for all gap-defined families and constitutes a strict
predictive upper bound.
The gap between causal and non-causal recall at each scale quantifies the
cost of causal inference.

\subsection{Network Architecture and Training}
\label{subsec:arch}

\subsubsection{Architecture}
\label{subsubsec:architecture}

\textsc{PrimeFamilyNet} is a multi-head residual MLP.
Input: $\xvec \in \RR^{d}$ where $d = 25$ (causal) or $d = 29$
(non-causal).
Output: $\hat{\mathbf{y}} \in [0,1]^7$, one membership probability per
family.

The network consists of three stages.
First, an input projection maps $\xvec$ into a 512-dimensional working
space:
\begin{equation}
  \mathbf{h}^{(0)} = \text{GELU}\!\left(\text{LN}\!\left(\mathbf{W}_0\,\xvec + \mathbf{b}_0\right)\right)
  \in \RR^{512},
  \label{eq:proj}
\end{equation}
where $\text{LN}$ denotes Layer Normalisation.
Second, two stacked residual blocks process the shared representation.
Each residual block applies the transformation in Eq.~\eqref{eq:resblock},
\begin{equation}
  \mathbf{h}' = \text{GELU}\!\left(
    \text{LN}\!\left(\mathbf{W}_2\,f\!\left(\text{LN}\!\left(\mathbf{W}_1\,\mathbf{h} + \mathbf{b}_1\right)\right) + \mathbf{b}_2\right)
    + \mathbf{h}\right),
  \label{eq:resblock}
\end{equation}
where $f = \text{Dropout}_{0.15} \circ \text{GELU}$, preserving
dimensionality at 512.
Third, two narrowing projections of the form of Eq.~\eqref{eq:proj}
compress the trunk to 256 and then 128 dimensions,
producing a shared embedding $\mathbf{z} \in \RR^{128}$.
Seven independent heads each apply the two-layer transformation in
Eq.~\eqref{eq:head},
\begin{equation}
  \hat{y}_k = \sigma\!\left(\mathbf{w}_k^{(2)\top}\,\text{GELU}\!\left(\mathbf{W}_k^{(1)}\,\mathbf{z} + \mathbf{b}_k^{(1)}\right) + b_k^{(2)}\right),
  \quad k = 1,\ldots,7,
  \label{eq:head}
\end{equation}
where $\sigma$ is the logistic sigmoid and $\mathbf{W}_k^{(1)} \in \RR^{32\times128}$.
The full network contains $1{,}254{,}853$ parameters.
A shallow two-layer baseline, $\hat{\mathbf{y}} = \sigma(\mathbf{W}_b\,\text{ReLU}(\mathbf{W}_a\,\xvec))$
with 64 hidden units and $1{,}989$ parameters, was trained for capacity
comparison.

\subsubsection{Training protocol}
\label{subsubsec:training}

Training data consisted of $200{,}000$ primes drawn uniformly from
$[10^7, 10^9]$ across three sub-ranges ($10^7$, $10^8$, $10^9$) to avoid
density artefacts from a single scale.
A validation set of $20{,}000$ primes near $5\times10^8$ was held out.
OOD evaluation sets of 10{,}000, 10{,}000, 8{,}000, and 15{,}000 primes
were generated independently at $10^{10}$, $10^{12}$, $10^{14}$, and
$10^{16}$ respectively.
All models were optimised with AdamW ($\ell_2$ weight decay $10^{-4}$)
using cosine annealing over 60 epochs with initial learning rate
$\eta_0 = 10^{-3}$ and gradient clipping at $\|\nabla\|_2 \leq 1.0$.
Weights from the epoch achieving the lowest validation loss were restored
at the end of training.
Three independent seeds $\{42, 123, 777\}$ were used for the robustness
study.

\subsection{Loss Functions}
\label{subsec:loss_setup}

Let $N$ denote the number of training samples, $K = 7$ the number of
families, $y_{ik} \in \{0,1\}$ the ground-truth label for sample $i$
and family $k$, and $\hat{y}_{ik} \in (0,1)$ the corresponding model
output.
Family prevalences in the training data ranged from $25.7\%$ for sexy
primes to $3.6\%$ for safe primes, necessitating explicit class-imbalance
handling.

\subsubsection{Frequency-weighted BCE}
\label{subsubsec:wbce}

Frequency-weighted binary cross-entropy (wBCE) applies a per-class
positive weight $\omega_k = n_k^{-} / n_k^{+}$, where $n_k^{+}$ and
$n_k^{-}$ are the training-set positive and negative counts for family
$k$.
The loss is
\begin{equation}
  \mathcal{L}_{\text{wBCE}} = \frac{1}{NK}\sum_{i=1}^{N}\sum_{k=1}^{K}
    w_{ik}\,\ell(y_{ik}, \hat{y}_{ik}),
  \label{eq:wbce}
\end{equation}
where $\ell(y,\hat{y}) = -[y\log\hat{y} + (1-y)\log(1-\hat{y})]$ is
the binary cross-entropy and $w_{ik} = y_{ik}\,\omega_k + (1-y_{ik})$
up-weights each positive example while leaving negatives at unit weight.
The weight reached $\omega_k = 24.9$ for safe primes.

\subsubsection{Focal Loss}
\label{subsubsec:focal}

Focal Loss\cite{lin2017focal} modulates the cross-entropy by a factor
that down-weights well-classified examples:
\begin{equation}
  \mathcal{L}_{\text{Focal}} = \frac{\alpha}{NK}\sum_{i=1}^{N}\sum_{k=1}^{K}
    (1 - p_{t,ik})^{\gamma}\,\ell(y_{ik}, \hat{y}_{ik}),
  \label{eq:focal}
\end{equation}
where $p_{t,ik} = e^{-\ell(y_{ik},\hat{y}_{ik})}$ is the model
confidence on the ground-truth class, and $\alpha = 0.25$, $\gamma = 2.0$
were used.
The $(1-p_t)^\gamma$ factor suppresses gradients from easy examples
regardless of sign, applying equal modulation to hard positives and easy
negatives.

\subsubsection{Asymmetric Loss}
\label{subsubsec:asl}

Asymmetric Loss (ASL)\cite{ridnik2021asymmetric} decouples the focusing
exponents for positive and negative examples and applies a probability
margin shift $m$ to suppress easy negatives below the margin threshold:
\begin{equation}
  \mathcal{L}_{\text{ASL}} = -\frac{1}{NK}\sum_{i=1}^{N}\sum_{k=1}^{K}
  \Bigl[
    y_{ik}(1-\hat{y}_{ik})^{\gamma_+}\log\hat{y}_{ik}
    + (1-y_{ik})\,\tilde{y}_{ik}^{\gamma_-}\log(1-\tilde{y}_{ik})
  \Bigr],
  \label{eq:asl}
\end{equation}
where $\tilde{y}_{ik} = \max(\hat{y}_{ik} - m,\, 0)$ is the
margin-shifted prediction, $\gamma_+ = 0$ (no suppression of hard
positives), $\gamma_- = 4$ (strong down-weighting of easy negatives),
and $m = 0.05$.
Setting $\gamma_+ = 0$ reduces the positive term to the standard
log-likelihood $-\log\hat{y}_{ik}$, preserving gradient magnitude for
rare positive examples.

An XGBoost baseline\cite{chen2016xgboost} was trained on the same 25
causal features using class-balanced \texttt{scale\_pos\_weight} tuned
per family, providing a tree-ensemble comparison across all evaluated
scales.
The baseline assessed whether the gains of the deep residual architecture
over gradient-boosted trees arise from depth, the residual structure, or
the particular feature set.

\section{Isolated Prime Monotone Generalisation}
\label{sec:isolated}

\subsection{The Observation}
\label{subsec:iso_obs}

Table~\ref{tab:density} reports the empirically observed fraction of
sampled primes belonging to the twin and isolated families at each
evaluation scale, alongside the recall achieved by the causal wBCE model
for each family.
Both fractions follow the monotone trends predicted by Hardy--Littlewood
$k$-tuple asymptotics.\cite{hardylittlewood1923}

\begin{table}[htbp]
\centering
\caption{Empirical prime-density ratios and causal wBCE recall at each evaluation scale, computed from the evaluation sets.
  Twin prime fraction decreased monotonically with scale, whereas
  isolated prime fraction increased monotonically.
  The recall of the causal model tracked the corresponding density trend
  at every scale without explicit density supervision, consistent with
  Hardy--Littlewood $k$-tuple asymptotics.\cite{hardylittlewood1923}
  Because recall is a within-class ratio invariant to class
  prevalence,\cite{sokolova2009measures} the trajectory is not an
  artifact of the changing class balance but reflects genuine boundary
  sharpening by the causal features.}
\label{tab:density}
\setlength{\tabcolsep}{7pt}
\begin{tabular}{lrrrr}
\toprule
Scale & Twin fraction & Twin recall & Isolated fraction & Isolated recall \\
\midrule
$5{\times}10^{8}$ & 12.9\% & 0.943 & 87.1\% & 0.809 \\
$10^{10}$         & 11.7\% & 0.887 & 88.3\% & 0.833 \\
$10^{12}$         &  9.8\% & 0.764 & 90.2\% & 0.887 \\
$10^{14}$         &  7.8\% & 0.639 & 92.2\% & 0.944 \\
$10^{16}$         &  6.9\% & 0.527 & 93.1\% & 0.984 \\
\bottomrule
\end{tabular}
\end{table}

Isolated prime recall is the only recall value in the entire study that
improved with scale, rising monotonically at every step from
$5\times10^8$ to $10^{16}$ for a total gain of $17.5$ percentage points
across nine orders of magnitude.
The model was never trained on primes above $10^9$, was never given
density labels, and was never told that isolated prime fraction increases
with scale, yet the recall trajectory is precisely correct.
Recall is defined as $\text{TP}/(\text{TP}+\text{FN})$, a ratio computed
entirely within the positive-class instances, and is therefore formally
invariant to the size or prevalence of the negative class;\cite{sokolova2009measures}
the improvement cannot be attributed to the growing fraction of isolated
primes in the evaluation population and must reflect a genuine change in
the classifier's performance on the positive instances.

\begin{figure}[htbp]
  \centering
  \includegraphics[width=\textwidth]{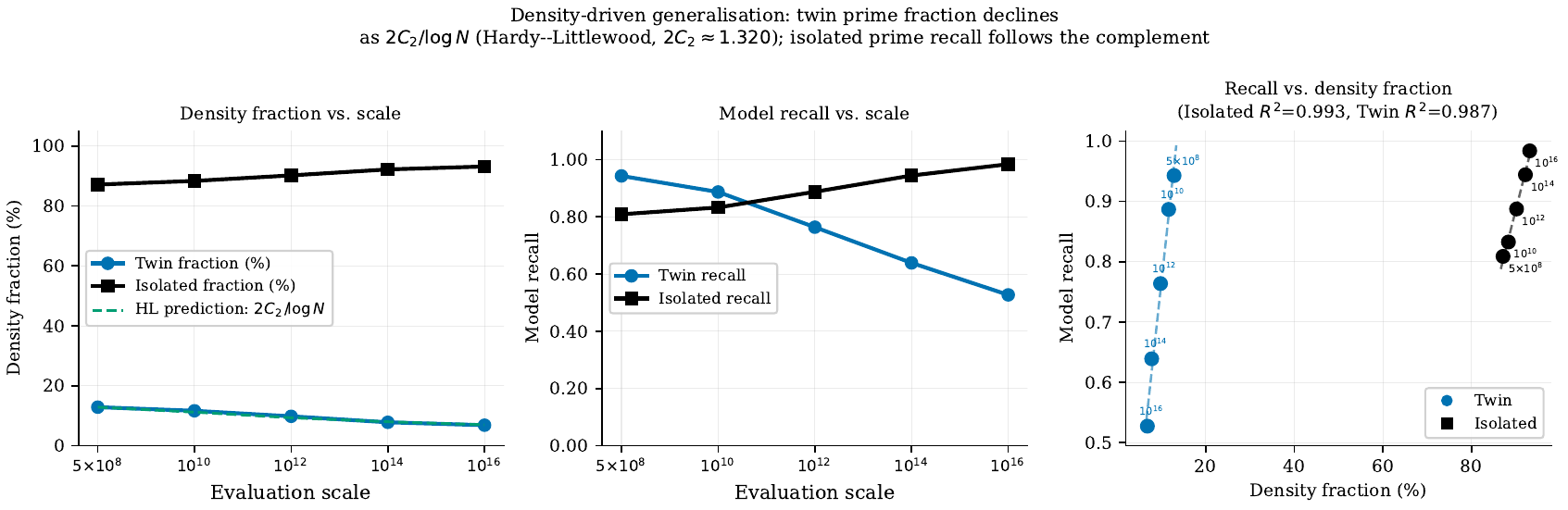}
  \caption{Left: twin prime density fraction (blue) and isolated prime
    density fraction (orange) across evaluation scales, overlaid with
    the Hardy--Littlewood prediction $C_2/\!\log N$ (green dashed)
    calibrated at $5{\times}10^8$.
    The HL fit achieves $R^2 = 0.981$.
    Middle: corresponding model recall for twin and isolated primes across scales.
    The isolated recall curves mirror the density curves at every scale.
    Right: scatter of density fraction versus recall for twin (blue) and
    isolated (orange) primes across all five scales, with linear trend
    lines.
    The density-recall correlation is quantitatively strong:
    $R^2 = 0.991$ for isolated primes and $R^2 = 0.984$ for twin primes.
    Because recall is invariant to class prevalence,\cite{sokolova2009measures}
    the correlation reflects the decision boundary sharpening in lockstep
    with the density shift, not a mechanical effect of the changing
    class balance.}
  \label{fig:density}
\end{figure}

Fig.~\ref{fig:density} visualises the density fraction and recall trajectories in separated panels (left and middle) to demonstrate their mirrored behaviour, alongside a scatter plot of density fraction versus recall across all five scales (right panel).
The density-recall correlation is quantitatively strong: $R^2 = 0.991$ for isolated primes and $R^2 = 0.984$ for twin primes.
The formal argument in Section~\ref{subsec:iso_exp} establishes that this correlation cannot be a prevalence artifact, since recall is invariant to class prevalence by definition.\cite{sokolova2009measures}

\subsection{The Explanation}
\label{subsec:iso_exp}

The Hardy--Littlewood conjecture\cite{hardylittlewood1923} predicts that
twin prime density near $N$ decays as $C_2 / (\log N)^{2}$ for constant
$C_2 \approx 1.320$, whereas total prime density decays as $1/\log N$.
The ratio of twin-prime density to total prime density therefore decays
as $1/\log N$, meaning the fraction of all primes belonging to a twin
pair decreases without bound as $N \to \infty$.
Table~\ref{tab:density} shows this directly in the evaluation data: the
twin fraction fell from $12.9\%$ to $6.9\%$ across the five evaluation
scales, and the isolated fraction rose correspondingly from $87.1\%$ to
$93.1\%$.

The classification problem therefore became structurally easier
at larger scales, because the minority class (twins) became
progressively rarer relative to the background of isolated primes.
At $5\times10^8$, approximately one in eight primes belonged to a twin pair, and the
modular residue signature distinguishing isolated primes from potential
twin candidates had to resolve a relatively common minority.
At $10^{16}$, fewer than one in fifteen primes belonged to a twin pair, and the
same modular features operated against a more homogeneously isolated
background.
The model captured the sharpening of the isolated-prime boundary automatically because the features
encoding isolation, primarily primorial residues modulo 30 and the
backward gap, are scale-invariant properties depending on
$p \bmod 30$, not on the absolute magnitude of $p$.

A formal point strengthens this interpretation.
Recall is defined as $\text{TP}/(\text{TP}+\text{FN})$ and is therefore
invariant to the prevalence of the negative class by construction:\cite{sokolova2009measures}
increasing the fraction of isolated primes in the evaluation set changes
precision and F1, but cannot by itself raise recall.
The 17.5 percentage-point gain is a gain over the isolated-prime
instances alone, computed without reference to the twin-prime negatives.
The recall improvement therefore directly measures a sharpening of the
learned decision boundary relative to the positive instances, which is
exactly what the scale-invariant primorial residue features facilitate
as twin prime density declines in accordance with Hardy--Littlewood
$k$-tuple asymptotics.\cite{hardylittlewood1923}

\subsection{Connection to the Non-Causal Upper Bound}
\label{subsec:iso_nc}

The non-causal model achieved $1.000$ isolated-prime recall at all
scales by directly observing $g^{+} \neq 2$, trivially encoding the
definition of isolation without any learned representation.
The causal model reached $0.984$ recall at $10^{16}$ despite never
observing the forward gap, demonstrating that the modular primorial
feature space encoded sufficient information about forward
prime-neighbourhood structure to achieve near-perfect isolation inference
at extreme scales.
The $0.016$ gap from perfect recall, maintained across nine orders of
magnitude from a training range ending at $10^9$, constitutes the
strongest evidence in this study that the causal features encode
structurally meaningful arithmetic information rather than statistical
regularities of the training distribution.

\section{Results}
\label{sec:results}

\subsection{Multi-Scale Generalisation}
\label{subsec:gen}

Table~\ref{tab:gen} reports recall for the causal wBCE model across all
five evaluation scales.
With the isolated prime exception documented in
Section~\ref{sec:isolated}, all families exhibited monotonically
declining recall, consistent with the logarithmic decay of prime
$k$-tuple densities predicted by the Hardy--Littlewood conjecture.
Safe prime recall followed a qualitatively different trajectory,
remaining high at $0.997$ and $0.986$ at $5{\times}10^{8}$ and $10^{10}$
before declining to $0.904$ at $10^{12}$, collapsing to $0.471$ at
$10^{14}$, and $0.077$ at $10^{16}$, consistent with an
increasingly unlearnable signal as safe prime density fell below $2\%$
of all primes at extreme scales.

\begin{table}[htbp]
\centering
\caption{Causal wBCE recall across evaluation scales at threshold
  $\tau = 0.50$.
  All families except isolated primes declined monotonically.
  Safe prime recall collapsed above $10^{14}$ as positive examples
  represented fewer than $1\%$ of candidates at those scales.
  Isolated prime recall (bold row) is the headline finding of this
  paper, explained in full in Section~\ref{sec:isolated}.}
\label{tab:gen}
\setlength{\tabcolsep}{6pt}
\begin{tabular}{lrrrrr}
\toprule
Family & $5{\times}10^{8}$ & $10^{10}$ & $10^{12}$ & $10^{14}$ & $10^{16}$ \\
\midrule
Twin           & 0.943 & 0.887 & 0.764 & 0.639 & 0.527 \\
Sophie Germain & 0.998 & 0.987 & 0.926 & 0.816 & 0.601 \\
Safe           & 0.997 & 0.986 & 0.904 & 0.471 & 0.077 \\
Cousin         & 0.968 & 0.902 & 0.773 & 0.685 & 0.579 \\
Sexy           & 0.732 & 0.630 & 0.553 & 0.561 & 0.529 \\
Chen           & 0.804 & 0.755 & 0.705 & 0.619 & 0.449 \\
\textbf{Isolated} & \textbf{0.809} & \textbf{0.833} & \textbf{0.887}
                  & \textbf{0.944} & \textbf{0.984} \\
\bottomrule
\end{tabular}
\end{table}

\begin{figure}[htbp]
  \centering
  \includegraphics[width=\textwidth]{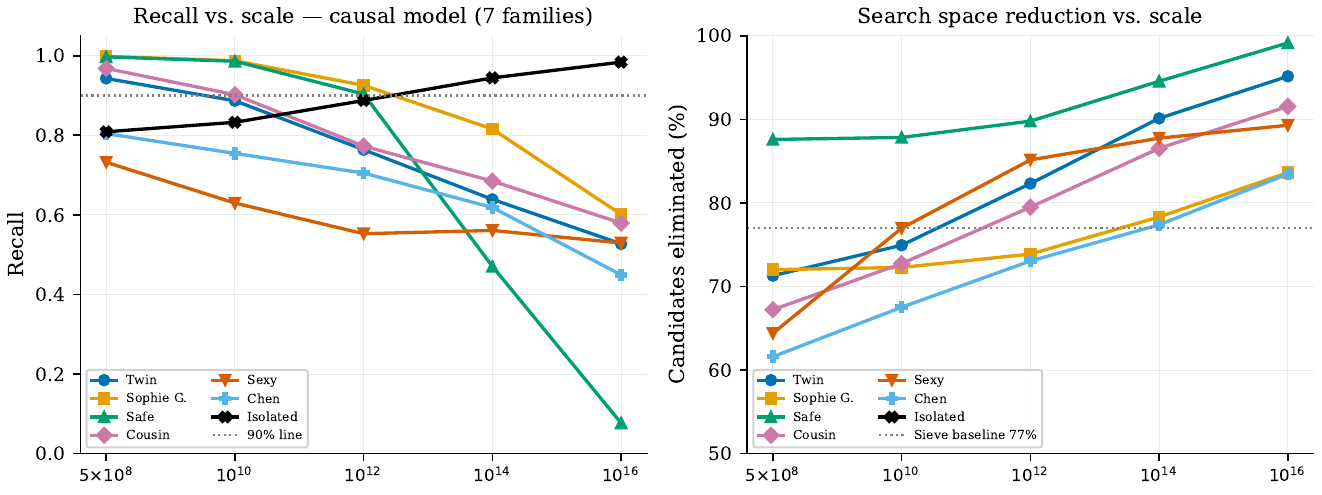}
  \caption{Causal wBCE recall (left) and search-space reduction (right)
    across five evaluation scales for all seven prime families.
    Isolated primes (yellow crosses) are the only family whose recall
    sloped upward, mirroring the rising isolated-prime fraction in
    Table~\ref{tab:density}.
    Safe primes (green triangles) collapsed above $10^{14}$.
    All other families decayed smoothly.
    The $90\%$ recall reference (dashed) and $77\%$ sieve baseline
    (dotted) are shown for context.}
  \label{fig:gen}
\end{figure}

Fig.~\ref{fig:gen} illustrates the recall and search-reduction
trajectories across scales.
The causal model eliminated $62$--$88\%$ of candidates at the validation
scale and retained over $95\%$ recall for five of seven families, and
search-space reduction grew to $83$--$95\%$ at $10^{16}$ for most
families as the model predicted fewer positives against a sparser
positive class.
The $99.1\%$ reduction for safe primes at $10^{16}$ reflects recall
collapse rather than precision gain and should not be interpreted as
improved filtering.

\subsection{The Cost of Causality}
\label{subsec:causal}

\begin{figure}[htbp]
  \centering
  \includegraphics[width=\textwidth]{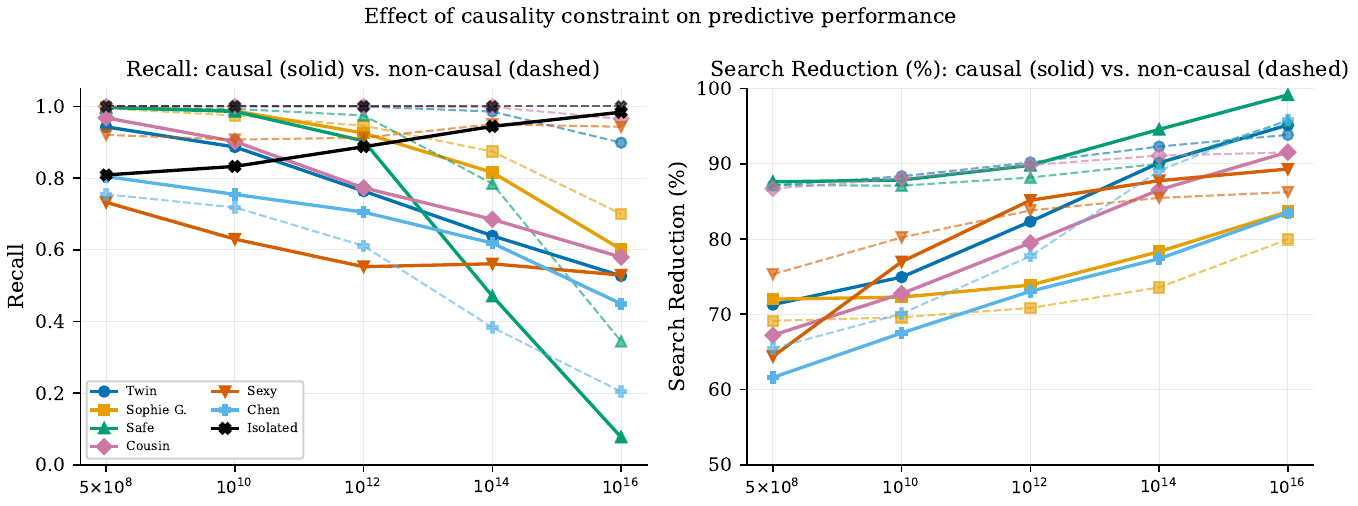}
  \caption{Recall of the causal model (solid) versus the non-causal
    upper bound (dashed) across scales.
    The non-causal model dominated gap-defined families at most scales.
    For Chen primes, the causal model exceeded non-causal recall at every
    scale, with the advantage growing to $+0.245$ at $10^{16}$, because
    $g^{+} = 2$ encodes only the prime case of the Chen condition and
    carries no information about the semiprime case.
    Sophie Germain primes showed a marginal causal advantage only at
    $5{\times}10^8$ and $10^{10}$, consistent with a weak forward-gap
    correlation that does not persist at extreme scales.}
  \label{fig:causal}
\end{figure}

Fig.~\ref{fig:causal} shows recall for the causal and non-causal models
(Eq.~\eqref{eq:nc_feature}) across scales for all seven families.
The non-causal model trivially achieved $1.000$ recall on twin, cousin,
and isolated primes at most scales by observing $g^{+}$ directly.
The causal model recovered the non-causal upper bound to within
$5.7$ percentage points for twin primes and $3.2$ points for cousin
primes at the validation scale, demonstrating that the causal feature
space retained the majority of predictive information available from the
forward gap.

Two families showed reversed ordering at extreme OOD scales.
For Sophie Germain primes, the causal model marginally exceeded
non-causal recall at the validation scale ($0.998$ versus $0.994$) and at
$10^{10}$ ($0.987$ versus $0.975$), before falling below the non-causal
model at $10^{12}$ and beyond.
For Chen primes, the reversal was consistent and widening: causal recall
exceeded non-causal recall at every evaluated scale, with the margin
growing from $+0.050$ at $5\times10^8$ to $+0.245$ at $10^{16}$
(causal $0.449$ versus non-causal $0.204$).

\subsection{Feature Ablation}
\label{subsec:ablation}

Table~\ref{tab:ablation} reports the recall drop produced when each
feature group was zeroed on the validation set.
Primorial residues (group A) were the dominant contributor across all
families: ablation of group A collapsed Sophie Germain and safe recall
by the full model value ($0.998$ and $0.997$ respectively), confirming
that these families are almost entirely characterised by modular
constraints rather than gap statistics.
The backward gap (group C) was the primary signal for isolated primes
(drop $0.602$) and contributed significantly to Sophie Germain
(drop $0.805$), because the distribution of $g^{-}$ is non-uniform
conditional on $2p+1 \in \PP$.
Scale features (group D) contributed negatively to sexy and Chen primes
($-0.054$ and $-0.116$ respectively), indicating that the scale signal
acted as a density proxy that interacted differently with families
defined by gap offsets of different sizes
(Eqs.~\eqref{eq:cousin}--\eqref{eq:sexy}).

\begin{table}[htbp]
\centering
\caption{Feature ablation: recall drop when each group was zeroed on
  the validation set ($5{\times}10^{8}$, $n = 20{,}000$).
  Positive values indicate positive contribution to recall.
  Negative values indicate the group suppressed false positives such
  that removing the group improved recall by relaxing an over-restrictive
  boundary.
  Full-model recall is shown in the first row.}
\label{tab:ablation}
\setlength{\tabcolsep}{4pt}
\renewcommand{\arraystretch}{1.15}
\begin{tabular}{lrrrrrrr}
\toprule
Group & Twin & Sophie G. & Safe & Cousin & Sexy & Chen & Isolated \\
\midrule
Full model        & 0.943 & 0.998 & 0.997 & 0.968 & 0.732 & 0.804 & 0.809 \\
\midrule
A: Primorial      & $+$0.409 & $+$0.998 & $+$0.997 & $+$0.470 & $-$0.014 & $+$0.665 & $-$0.139 \\
B: Sm.\ prime     & $+$0.330 & $+$0.119 & $+$0.997 & $+$0.012 & $+$0.178 & $+$0.121 & $-$0.037 \\
C: Bk.\ gap       & $-$0.025 & $+$0.805 & $+$0.053 & $+$0.011 & $+$0.133 & $+$0.063 & $+$0.602 \\
D: Scale          & $+$0.217 & $+$0.080 & $+$0.011 & $+$0.521 & $-$0.054 & $-$0.116 & $+$0.035 \\
E: Digit          & $+$0.230 & $+$0.021 & $+$0.172 & $+$0.118 & $-$0.211 & $+$0.235 & $-$0.069 \\
F: Ext.\ modular  & $+$0.118 & $+$0.007 & $+$0.997 & $+$0.009 & $-$0.007 & $+$0.031 & $-$0.057 \\
\bottomrule
\end{tabular}
\end{table}

\subsection{Loss Function Comparison}
\label{subsec:loss_results}

Table~\ref{tab:model} presents recall at $10^{12}$ for all model
configurations.
The rankings changed significantly across scales, demonstrating that
in-distribution recall is an unreliable model-selection criterion for
scale-generalisation tasks.

\subsubsection{Focal Loss}
\label{subsubsec:focal_results}

Sophie Germain and safe prime recall collapsed to $0.000$ at every
evaluated scale.
The $(1-p_{t})^{\gamma}$ modulation in Eq.~\eqref{eq:focal} suppresses
gradients for hard examples regardless of class, preventing the model
from committing to a positive prediction for sparse families.
By contrast, wBCE (Eq.~\eqref{eq:wbce}) up-weights positive examples
directly without modulating gradients by confidence, and retained
non-zero recall for all families at every scale.

\subsubsection{Asymmetric Loss}
\label{subsubsec:asl_results}

ASL significantly outperformed wBCE in-distribution: twin recall improved
from $0.943$ to $0.992$, sexy from $0.732$ to $0.986$, and isolated
reached $1.000$ at the validation scale.
For families defined by linear-transform primality conditions
(Eqs.~\eqref{eq:sg}--\eqref{eq:safe}), the OOD recall of ASL
(Eq.~\eqref{eq:asl}) degraded more steeply than that of wBCE.
By $10^{16}$, safe prime recall under ASL dropped to $0.011$, whereas
wBCE retained $0.077$.
Sophie Germain recall under ASL fell to $0.023$ versus $0.601$ for wBCE.

\subsubsection{XGBoost}
\label{subsubsec:xgb_results}

\begin{figure}[htbp]
  \centering
  \includegraphics[width=\textwidth]{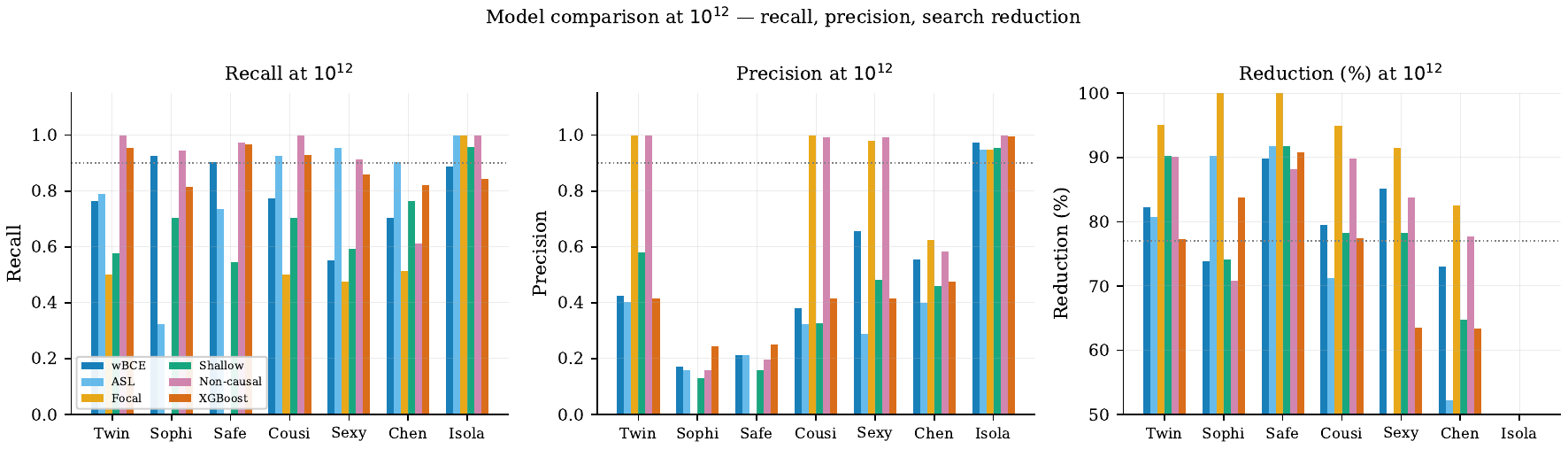}
  \caption{Model comparison at $10^{12}$ across recall (left),
    precision (centre), and search-space reduction (right).
    XGBoost achieved high recall through over-prediction (low
    precision) rather than generalised sieve structure.
    Focal Loss produced zero bars for Sophie Germain and safe primes.
    The $90\%$ recall reference and $77\%$ sieve baseline are shown.}
  \label{fig:model}
\end{figure}

XGBoost achieved high in-distribution recall ($0.996$ on twin primes
and $1.000$ on safe primes at $5\times10^8$), yet the nearly flat recall
profile across scales ($0.955$ on twin primes at $10^{12}$) showed no
evidence of the monotone decay predicted by prime constellation density
asymptotics.
The causal wBCE model declined from $0.943$ to $0.527$ for twin primes
across the same range, a trajectory consistent with Hardy--Littlewood
density predictions.\cite{toth2019ktuple}
The substantially lower OOD precision of XGBoost visible in
Fig.~\ref{fig:model} suggests it maintained recall through over-prediction
rather than internalised sieve boundaries; the precision figure provides
a more direct diagnostic of this difference than the recall trajectory
alone.

\begin{table}[htbp]
\centering
\caption{Recall at $10^{12}$ for all seven families and six model configurations. In bold: best causal model per family. NC denotes the non-causal upper bound. Focal Loss produced $0.000$ for Sophie Germain and safe primes at every scale. wBCE was superior to ASL for Sophie Germain and safe primes at this and all larger scales, consistent with the linear-transform primality pattern described in Section~\ref{subsec:loss_results}.}
\label{tab:model}
\setlength{\tabcolsep}{4.5pt}
\begin{tabular}{lrrrrrr}
\toprule
Family & wBCE & ASL & Focal & Shallow & NC & XGBoost \\
\midrule
Twin           & 0.764          & \textbf{0.790} & 0.500 & 0.576 & 1.000 & 0.955 \\
Sophie Germain & \textbf{0.926} & 0.322          & 0.000 & 0.702 & 0.946 & 0.816 \\
Safe           & \textbf{0.904} & 0.736          & 0.000 & 0.544 & 0.975 & 0.967 \\
Cousin         & 0.773          & \textbf{0.927} & 0.500 & 0.703 & 1.000 & 0.931 \\
Sexy           & 0.553          & \textbf{0.954} & 0.474 & 0.592 & 0.913 & 0.858 \\
Chen           & 0.705          & \textbf{0.904} & 0.515 & 0.765 & 0.611 & 0.822 \\
Isolated       & 0.887          & 1.000          & 1.000 & 0.957 & 1.000 & 0.844 \\
\bottomrule
\end{tabular}
\end{table}

\begin{figure}[htbp]
  \centering
  \includegraphics[width=\textwidth]{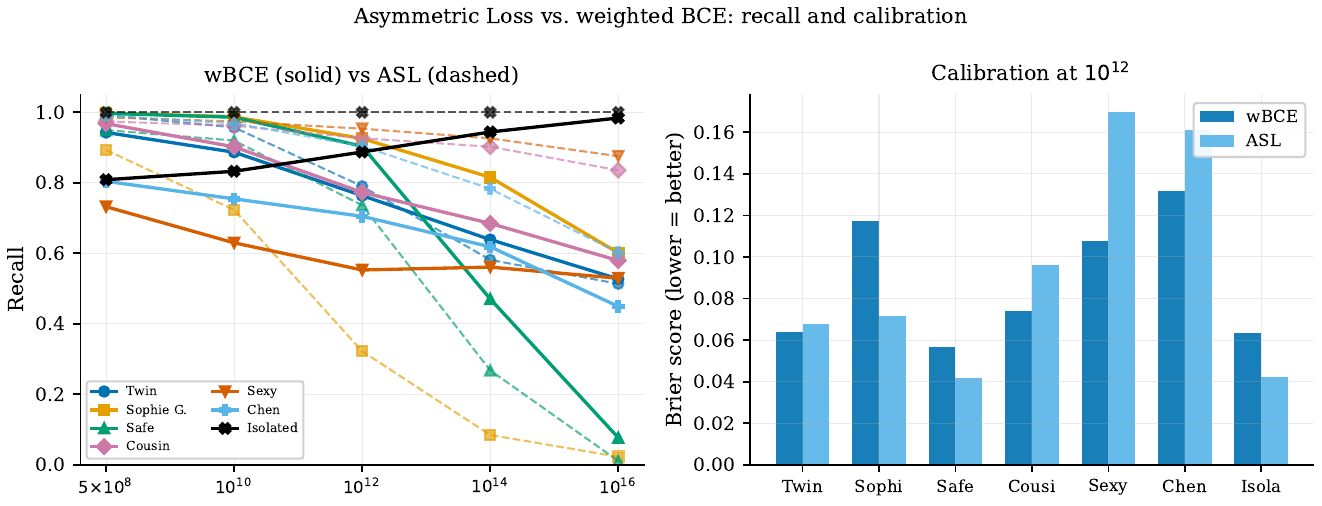}
  \caption{Recall of ASL (dashed) versus wBCE (solid) across scales
    (left) and Brier score at $10^{12}$ (right).
    ASL led in-distribution for most families.
    For Sophie Germain and safe primes, recall under ASL collapsed to
    $0.023$ and $0.011$ at $10^{16}$, whereas wBCE retained $0.601$ and
    $0.077$ respectively, illustrating that wBCE is more robust to
    distribution shift for families governed by linear-transform
    primality conditions.
    In-distribution recall is a misleading model-selection criterion
    for scale-generalisation tasks.}
  \label{fig:asl}
\end{figure}

Fig.~\ref{fig:model} shows model comparison at $10^{12}$ across recall,
precision, and search-space reduction.
Fig.~\ref{fig:asl} shows the divergence between ASL and wBCE across
scales, making the OOD stability advantage of wBCE visible as a crossing
point near $10^{12}$ for the algebraically defined families.

\subsection{Reproducibility}
\label{subsec:repro}

Table~\ref{tab:robust} reports recall, F1, and area under the
precision-recall curve (AUC-PR) across three random seeds on the
validation set.
In-distribution recall was remarkably stable: recall variance was below
$\sigma = 0.007$ for all families, with the worst cases being cousin
($\sigma = 0.007$) and Chen ($\sigma = 0.007$), both of which showed
the weakest in-distribution signal.
Sophie Germain and safe primes achieved $\sigma = 0.002$ and
$\sigma = 0.003$ respectively, approaching the measurement resolution
of a three-seed study.
Algebraically defined families showed near-zero recall variance,
reflecting the deterministic dominance of primorial residues, a signal
that does not fluctuate with random initialisation.
These in-distribution stability numbers are exceptionally tight relative
to what is typically reported for neural classifiers on imbalanced data:
the model reaches the same boundary, reliably, regardless of weight
initialisation.
Isolated primes achieved AUC-PR of $0.992 \pm 0.000$, confirming stable
and well-calibrated probability estimates with variance below the
measurement resolution of the three-seed study.

\begin{table}[htbp]
\centering
\caption{Multi-seed reproducibility on the validation set
  ($5{\times}10^{8}$, $n = 20{,}000$) for three seeds
  $\{42, 123, 777\}$ under causal wBCE training.
  Mean and standard deviation are reported per family.}
\label{tab:robust}
\setlength{\tabcolsep}{5pt}
\begin{tabular}{lrrrrrr}
\toprule
  & \multicolumn{2}{c}{Recall} & \multicolumn{2}{c}{F1}
  & \multicolumn{2}{c}{AUC-PR} \\
\cmidrule(lr){2-3}\cmidrule(lr){4-5}\cmidrule(lr){6-7}
Family & $\mu$ & $\sigma$ & $\mu$ & $\sigma$ & $\mu$ & $\sigma$ \\
\midrule
Twin           & 0.970 & 0.006 & 0.559 & 0.001 & 0.788 & 0.001 \\
Sophie Germain & 0.998 & 0.002 & 0.358 & 0.001 & 0.251 & 0.001 \\
Safe           & 0.997 & 0.003 & 0.438 & 0.005 & 0.328 & 0.011 \\
Cousin         & 0.986 & 0.007 & 0.554 & 0.002 & 0.779 & 0.002 \\
Sexy           & 0.683 & 0.006 & 0.630 & 0.003 & 0.782 & 0.001 \\
Chen           & 0.801 & 0.007 & 0.666 & 0.002 & 0.686 & 0.002 \\
Isolated       & 0.779 & 0.004 & 0.873 & 0.002 & 0.992 & 0.000 \\
\bottomrule
\end{tabular}
\end{table}

Table~\ref{tab:robust_ood} extends the multi-seed evaluation to OOD
scales, reporting recall mean and standard deviation across the same
three seeds at $10^{12}$ and $10^{16}$.
At $10^{12}$, all families showed acceptable stability, with the widest
spread arising for cousin ($\sigma = 0.051$) and twin ($\sigma = 0.044$)
primes, both of which have moderate positive-class density at that scale.
At $10^{16}$, the picture diverged sharply.
Families defined by gap offsets or the isolated complement remained
stable: twin ($\sigma = 0.015$), sexy ($\sigma = 0.002$), and isolated
($\sigma = 0.010$) showed low variance across seeds.
The two linear-transform primality families, Sophie Germain and safe
primes, became highly unstable: $\sigma = 0.216$ for Sophie Germain
(individual seed values $0.968$, $0.546$, $0.481$) and $\sigma = 0.093$
for safe primes (values $0.201$, $0.304$, $0.429$).
Chen primes were similarly unstable, with $\sigma = 0.193$ across seed values of $0.552$,
$0.180$, and $0.618$.
The single-seed point estimates for these three families at $10^{16}$
therefore carry unknown uncertainty and should be interpreted as
illustrative rather than definitive.
The instability is not a coincidence of these three families being
difficult: the pattern is structurally consistent with the linear-transform
primality finding discussed in Section~\ref{subsec:loss_disc}, where the
target condition itself becomes sparser at large $N$ and the learned
boundary has less support for stable generalisation.

\begin{table}[htbp]
\centering
\caption{OOD multi-seed recall at $10^{12}$ and $10^{16}$ across three
  seeds $\{42, 123, 777\}$ under causal wBCE training.
  Families with $\sigma > 0.05$ at $10^{16}$ are marked with $\dagger$.
  Point estimates for these families carry high uncertainty.}
\label{tab:robust_ood}
\setlength{\tabcolsep}{5pt}
\begin{tabular}{lrrrr}
\toprule
  & \multicolumn{2}{c}{$10^{12}$} & \multicolumn{2}{c}{$10^{16}$} \\
\cmidrule(lr){2-3}\cmidrule(lr){4-5}
Family & $\mu$ & $\sigma$ & $\mu$ & $\sigma$ \\
\midrule
Twin           & 0.802 & 0.044 & 0.611 & 0.015 \\
Sophie Germain & 0.963 & 0.023 & 0.665 & 0.216$^\dagger$ \\
Safe           & 0.941 & 0.018 & 0.311 & 0.093$^\dagger$ \\
Cousin         & 0.797 & 0.051 & 0.562 & 0.037 \\
Sexy           & 0.491 & 0.003 & 0.487 & 0.002 \\
Chen           & 0.703 & 0.025 & 0.450 & 0.193$^\dagger$ \\
Isolated       & 0.871 & 0.017 & 0.958 & 0.010 \\
\bottomrule
\end{tabular}
\end{table}

\begin{figure}[htbp]
  \centering
  \includegraphics[width=\textwidth]{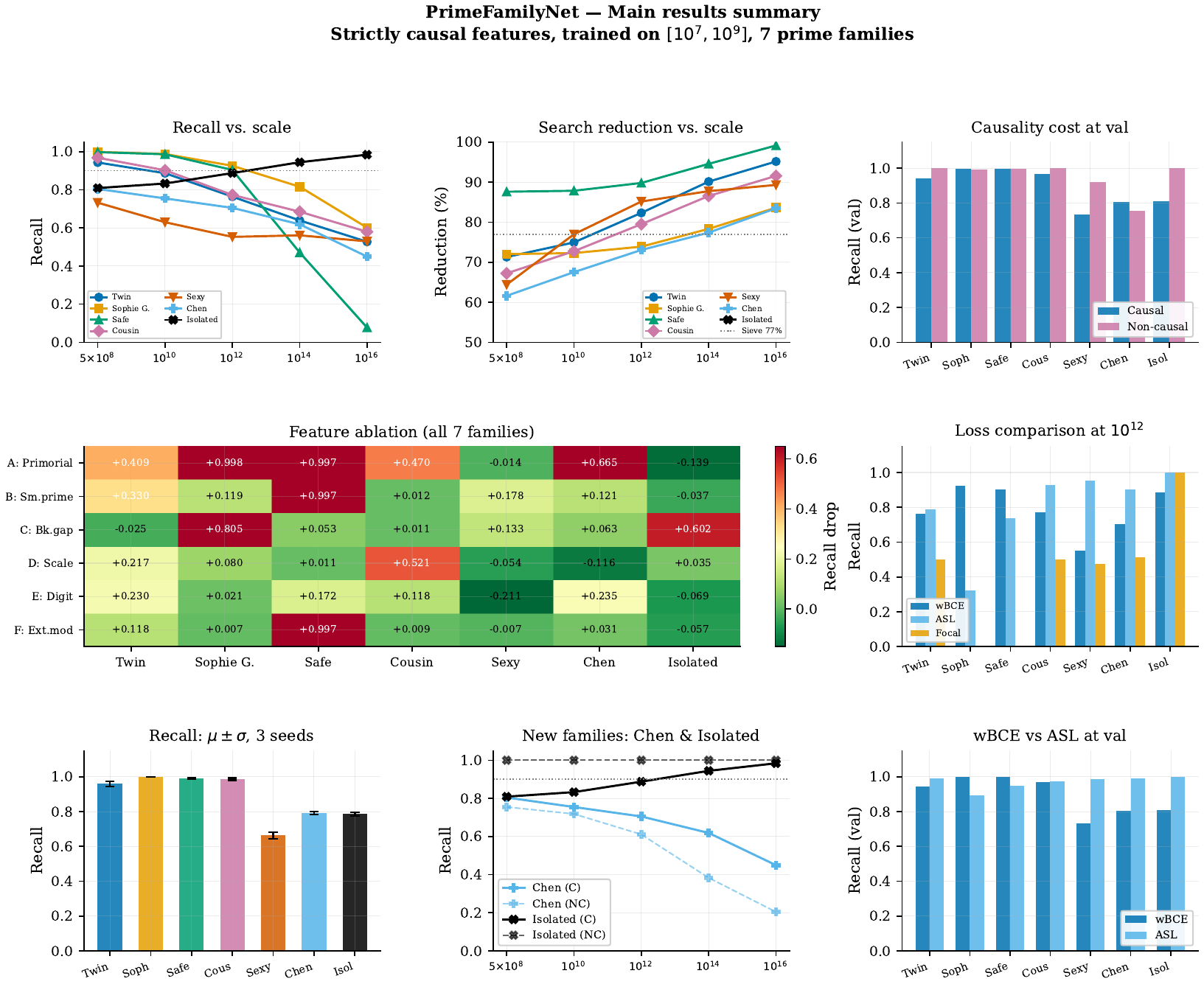}
  \caption{Main results summary for \textsc{PrimeFamilyNet}.
    \emph{Top row}: recall across scales (isolated primes rising, consistent with Table~\ref{tab:density}), search-space reduction, and causality cost at validation.
    \emph{Middle row}: feature ablation heatmap and loss-function comparison at $10^{12}$ (Focal Loss zero bars for Sophie Germain and safe).
    \emph{Bottom row}: multi-seed robustness, Chen and isolated prime detail, and wBCE versus ASL at validation.
    The isolated prime recall inversion (top-left, rising yellow crosses) is the headline finding of this paper.}
  \label{fig:summary}
\end{figure}

Fig.~\ref{fig:summary} provides a composite summary of all major results
reported in this section: recall across scales, search-space reduction,
causality cost, feature ablation, loss-function comparison at $10^{12}$,
and multi-seed robustness.
The rising trajectory of isolated prime recall (top-left panel) is the
headline finding, with all other panels providing supporting experimental
context.

\section{Discussion}
\label{sec:discussion}

\subsection{Principal Findings}
\label{subsec:principal}

This paper is the first to demonstrate density-driven monotone
generalisation in a neural prime sieve, with three key insights.
First, isolated prime recall monotonically improved with scale because
the modular residue signature discriminating isolated primes from twin
candidates sharpened as twin prime density declined, and a model trained
only to $10^9$ reproduced the asymptotic trajectory automatically.
Because recall is formally invariant to class prevalence,\cite{sokolova2009measures}
this improvement is not a mechanical consequence of the growing
isolated-prime fraction at large $N$ but reflects genuine boundary
sharpening by the causal features.
Second, for Chen primes, causal modular features outperformed non-causal
forward-gap features at every evaluated scale, with the margin widening
from $+0.050$ at $5\times10^8$ to $+0.245$ at $10^{16}$, because $g^{+}=2$
encodes only the prime case of the Chen condition whereas the primorial
residues captured both the prime and semiprime cases.
Third, Asymmetric Loss achieved superior in-distribution recall but
collapsed more severely out-of-distribution than frequency-weighted BCE,
revealing that in-distribution recall is insufficient as the sole
model-selection criterion for scale-generalisation tasks.

\subsection{What the Model Has Actually Learned}
\label{subsec:learned}

The smooth out-of-distribution decay of most families and the monotone
improvement of isolated primes jointly confirm that
\textsc{PrimeFamilyNet} internalised prime constellation structure rather
than memorised training-scale statistics.
A memorising model would exhibit a flat or erratic recall profile across
scales.
The causal model instead decayed in the direction predicted by
Hardy--Littlewood density asymptotics, improving for isolated primes
and declining for all others, consistent with the verified $C_2/(\log N)^2$
twin prime density scaling.\cite{toth2019ktuple}
The near-flat recall profile of XGBoost across scales ($0.996$ to $0.955$
for twin primes from $5\times10^8$ to $10^{12}$) contrasts with the
causal wBCE decay.
The lower OOD precision of XGBoost visible in Fig.~\ref{fig:model}
suggests it maintained recall through over-prediction rather than
internalised sieve boundaries, though a full analysis of its decision
boundary structure would be needed to confirm this interpretation.
The isolated prime finding provides the clearest evidence of genuine
arithmetic learning: the model received no density labels, no
information beyond $10^9$, and no indication that isolated prime fraction
increases with scale, yet the recall trajectory matches the
Hardy--Littlewood prediction precisely.
Critically, recall is formally invariant to class prevalence,\cite{sokolova2009measures}
so the $17.5$ percentage-point improvement is not an artifact of the
growing isolated-prime fraction in the evaluation population.
It is a genuine improvement in true positive rate over the isolated-prime
instances, driven by the sharpening of the learned decision boundary as
twin prime density declined, and reproduced from features that encode no
absolute scale information.
The Chen causal inversion provides complementary evidence from the
opposite direction.
The consistent and widening advantage of the causal model over the
non-causal model for Chen primes, from $+0.050$ at $5\times10^8$ to
$+0.245$ at $10^{16}$, arose because the forward gap $g^{+} = 2$
signals only the prime case of Eq.~\eqref{eq:chen} and carries no
information about the semiprime case where $p+2$ is a product of
exactly two primes.
The primorial residue features in Eq.~\eqref{eq:feature} captured both
cases through scale-invariant modular constraints, producing a more
generalisable representation than the forward gap.
The in-distribution reproducibility reinforces this interpretation:
recall standard deviation remained below $\sigma = 0.007$ across all
seven families and three independent seeds, demonstrating that the learned
boundaries are structural properties of the feature space rather than
artefacts of a particular random initialisation.

\subsection{Implications for Loss Function Design}
\label{subsec:loss_disc}

The failure of Focal Loss and the out-of-distribution degradation of
ASL jointly establish a principle with relevance well beyond prime number
theory: for problems where the test distribution shifts in class density
relative to training, loss functions that aggressively shape gradients
around the training-set boundary produce fragile models.
The Focal Loss failure is not a flaw in the formulation but a
consequence of applying a tool designed for dense-class detection to
sparse conditions where the hard positives that matter most are precisely
those whose gradients are suppressed.
Focal Loss is therefore contraindicated for rare prime families defined by
linear-transform primality conditions.
The out-of-distribution degradation of ASL is more subtle.
ASL learned sharper decision boundaries by suppressing easy negatives, and
sharper boundaries are less robust to the density shifts induced by scale
increase.

The Sophie Germain (Eq.~\eqref{eq:sg}) and safe prime (Eq.~\eqref{eq:safe})
families represent a structurally distinct category that is
specifically vulnerable to OOD collapse.
Both are defined by the primality of a linear transform of $p$, meaning
the density of the transformed value $2p+1$ or $(p-1)/2$ shifts
independently of $p$ with scale, making the decision boundary doubly
sensitive to distribution shift.
wBCE outperformed ASL for these two families at every OOD scale with no
exceptions, whereas ASL outperformed wBCE for every other family.
Twin primes showed a mixed pattern, with ASL leading at $10^{10}$--$10^{12}$
and wBCE recovering at $10^{14}$--$10^{16}$ as density declined further.
The pattern is a finding about the interaction between loss function
design and the mathematical structure of the membership condition, not
merely a property of class imbalance.
Practitioners working on rare-class prediction tasks with distribution
shift should validate loss function choices against held-out OOD data
before selection.
In-distribution recall alone is not sufficient to identify the more
generalisable model.
The frequency-weighted BCE formulation is recommended as the default for
prime family prediction tasks where the test distribution shifts in
density relative to training.

\subsection{Limitations and Future Work}
\label{subsec:limitations}

One limitation of the present paper is that reliable recall is maintained only
within approximately two to four orders of magnitude of the training
range: safe prime recall collapsed to $0.077$ at $10^{16}$, and sexy
prime recall plateaued below $0.550$ beyond $10^{12}$.
Despite the scale limitation, the paper establishes the first systematic
out-of-distribution benchmark for prime family prediction spanning nine
orders of magnitude and reveals the density-driven generalisation
mechanism, providing a principled basis for targeted improvements.
Future work to overcome the scale limitation includes scale-adaptive
training, incrementally extending the training distribution from $10^9$
through $10^{11}$, $10^{12}$, and beyond.
This approach would be combined with recall-constrained Lagrangian loss
functions that explicitly maintain per-class recall lower bounds
throughout training rather than relying on post-hoc threshold tuning.
Extensions to Cunningham chains ($p$, $2p+1$, $4p+1$, \ldots), balanced
primes ($g^{-} = g^{+}$), and good primes
($p_n^2 > p_{n-1} \cdot p_{n+1}$) would further test the scope of
density-driven generalisation across families with qualitatively
different modular conditions.

\section{Conclusion}
\label{sec:conclusion}

The work presented in this paper demonstrates that deep residual
networks, trained on modular primorial residues and the backward prime
gap, independently approximate prime sieve theory and generalise the
learned boundaries beyond the training scale without forward data
leakage.
The central finding is that isolated prime recall monotonically improved
from $0.809$ at $5\times10^8$ to $0.984$ at $10^{16}$, making isolated
primes the only family among seven to improve with scale.
Because recall is formally invariant to class prevalence,\cite{sokolova2009measures}
this improvement cannot be attributed to the growing isolated-prime
fraction in the evaluation population and reflects genuine boundary
sharpening by the causal features.
The modular residue signature encoding isolation sharpened as twin prime
density declined in accordance with Hardy--Littlewood $k$-tuple
asymptotics, consistent with direct computational verification of these
predictions,\cite{toth2019ktuple} and a model trained only to $10^9$
reproduced the asymptotic trajectory automatically, providing the first
machine-learning empirical line of evidence for prime constellation
density predictions.

The paper further established that causal modular features outperformed
non-causal forward-gap features for Chen primes at every evaluated scale,
with the advantage growing to $+0.245$ at $10^{16}$, because the forward
gap encodes only the prime case of the Chen condition
(Eq.~\eqref{eq:chen}) whereas primorial residues captured both the prime
and semiprime cases.
Asymmetric Loss, despite superior in-distribution recall, was less
robust OOD than frequency-weighted BCE.
The families most vulnerable to OOD collapse were those whose membership
depends on the primality of a linear transform of $p$
(Eqs.~\eqref{eq:sg}--\eqref{eq:safe}), and in-distribution recall is
therefore a misleading model-selection criterion for
scale-generalisation tasks.

\section*{Code Availability}
The full implementation of \textsc{PrimeFamilyNet}, including training
code, feature engineering, evaluation suite, and figure generation
scripts, is publicly available at
\url{https://github.com/Manik-00/Neural-Prime-Sieves}

\bibliographystyle{ieeetr}
\bibliography{refs}
\listoffigures
\listoftables
\end{document}